\documentclass[a4paper,12pt]{article}
\usepackage[cmex10]{amsmath}
\usepackage{latexsym,amssymb}
\usepackage[mathscr]{eucal}
\usepackage{amsxtra,amscd,amsthm,upref,layout,color}
\usepackage{calc}
\usepackage[final]{graphicx}
\def\refup#1{{$^{#1}$}}\def\refeq#1{(\ref{#1})}
\def\p{^{^{\prime}}}\def\pp{^{^{\prime\prime}}}

\def\min{{\rm {min}}}

\def\oo{{\leavevmode\setbox0=\hbox{h}\dimen0=\ht0 \advance\dimen0
by-1ex\rlap{\raise0.47\dimen0\hbox{\char'27}}o}}

\def\begeq{\begin{equation}}
\def\endeq{\end{equation}}
\def\begdis{\begin{displaymath}}
\def\enddis{\end{displaymath}}

\def\cA{{\cal A}}

\def\cS{{\cal S}}

\def\gotC{{\mathfrak{C}}}\def\gotD{{\mathfrak{D}}}

\def\gotG{{\mathfrak{G}}}\def\gotGdd{{\gotG}^{[2]}}

\def\S1{{S$_1$}}

\def\ie{{\em i.e.}}
\def\etal{{\em et al.}}

\begin{document}
\title{Algebraic Approximations of a Polyhedron Correlation Function 
Stemming from its Chord Length Distribution}
\author{ 
{{ 
 Salvino Ciccariello$^{a,b,*}$
}}\\
%
  \begin{minipage}[t]{0.9\textwidth}
   \begin{flushleft}
\setlength{\baselineskip}{12pt}
{\slshape  {\footnotesize{
$^{a}$ Universit\`{a} di Padova, Dipartimento di Fisica {\em {G. Galilei}}, 
 Via Marzolo 8, I-35131 Padova, Italy, 
and \\
$^{b}$Universit\`{a}  {\em {C\`{a} Foscari}}, Department of Molecular 
Sciences and Nanosystems, 
Via Torino 155/B, I-30172 Venezia, Italy
}}}\\
 \footnotesize{$^*$salvino.ciccariello@unipd.it}
\end{flushleft}
\end{minipage}
}      
%
\date{\today  }   
\maketitle                        
\centerline{\sl To be published in Acta Cryst. A{\bf77},..,(2021),} 
\centerline{\sl https://doi.org/10.1107/S2053273320014229} 
\begin{abstract} 
\noindent An  algebraic approximation, of order $K$,  of a polyhedron  
correlation function (CF)  can be obtained from $\gamma\pp(r)$, its 
chord-length distribution (CLD),  considering first, within the subinterval 
$[D_{i-1},\, D_i]$  of the full 
range of distances, a polynomial in the two variables $(r-D_{i-1})^{1/2}$ and 
$(D_{i}-r)^{1/2}$ such that its expansions around $r=D_{i-1}$ and $r=D_i$ 
simultaneously coincide with left and the right expansions of  $\gamma\pp(r)$ 
around $D_{i-1}$ and $D_i$ up to the terms $O\big(r-D_{i-1}\big)^{K/2}$ 
and $O\big(D_i-r\big)^{K/2}$, respectively. Then, for each $i$, one integrates twice  
the  polynomial  and determines the integration constants matching the resulting 
integrals at the common end points.  The 3D Fourier transform of the resulting 
algebraic CF approximation correctly reproduces, at large $q$s,  the asymptotic 
behaviour of the exact form  factor  up to the term $O(q^{-(K/2+4)})$. For 
illustration, the procedure is applied to  the cube, the tetrahedron and the 
octahedron.  \\      

{\em Synopsis} We report a procedure for obtaining an algebraic approximation of the correlation
function of a polyhedron starting from its known chord length distribution. 
    
\end{abstract}
\vfill

\eject
\noindent{\bf{1 -  Introduction } } \\
\indent Nowadays, the correlations functions (CF) [$\gamma(r)$] of the 
first three Platonic solids are explicitely known.  
In fact, Goodisman (1980) directly obtained the cube CF by 
evaluating the angular average of the overlapping volume while Ciccariello, 
starting from the general integral expression of the CF's 2nd order  derivative 
(Ciccariello \etal, 1981), first worked out the explicit expressions of the CLDs 
of the tetrahedron (2005a) and the octahedron (2014a),  and later succeeded in 
integrating these CLDs to get the corresponding CFs (2014b).   The expressions of 
these CFs are not simple and take different forms within each of  
the $M$ subintervals  of the full range of distances 
$[0,\, D_M]$.   \{For $i=1,\ldots,M$. the  $i$th subinterval is defined as 
$[D_{i-1},\,D_i]$ (with $D_0\equiv 0$), its length is denoted by 
$\Delta_{i}\equiv (D_i-D_{i-1})$ while the $D_i$s are some of the distance 
values between  vertices, 
between vertices and  sides, and  between  the sides of the given polyhedron.\}  
However, within each subinterval,  the CFs are analytic functions of $r$ and their 
structure  always is a  sum of rational functions  and of inverse trigonometric 
functions, the arguments of which  also are rational functions.  The last  
functions have the form $R(r,y)\equiv P(r,y)/Q(r,y)$ where $P$ and $Q$ 
are polynomials and $y$ denotesthe square root of a 2nd degree polynomial of $r$.   
Hence the 
important property:  the derivatives of this kind of functions, whatever their order,  
are functions of the same kind. Quite recently, Ciccariello (2020a,b) has shown that 
the mentioned mathematical structure also applies to the CLD of any bounded 
polyhedron, whatever its shape. However, the related  CFs are, as yet,  not 
eplicitly known owing to the difficulty of twice integrating  the CLDs in a closed 
analytic form. \\   
\indent In his report on the ms. of Ciccariello (2020b) paper one of the 
referees raised an important question, namely: whether it is possible to get an approximate algebraic 
expression of the CF stemming from the reported CLD. 
In this short note we present a procedure that achieves this aim through 
the following steps. Consider for definitenss the $i$th subinterval. Since we know
the analytic form of the CLD inside each subinterval, we also know its right and 
left expansions respectively around the two end points $D_{i-1}$ and $D_i$ of 
the considered interval. The truncation of the two expansions yield two 
algebraic expressions which respectively approximate the CLD around $D_{i-1}$ 
and $D_i$.  The difficulty now is that of devising a single algebraic function 
which simultaneously almost coincides with the truncated right expansion as 
$r\to D_{i-1}$ 
and with the left truncated expansion as $r\to D_i$. This problem is solved in the 
following section. The resulting function yields an algebraic approximation of 
the CLD within the full $i$th interval and its accuracy generally depends on the 
truncation order $K$. Integrating twice the resulting function we obtain an algebraic 
approximation of the CF within the same subinterval once we have 
determined the two integration constants. This determination is achieved by 
matching the CF approximations relevant to, say, the $i$th and the $(i+1)$th 
subintervals at the common end-point $D_i$. (The details of the procedure are 
reported in the second part of section 2.) Then, the sought for algebraic 
approximation of the full CF results from the combination of all these subinterval 
approximations. Section 3 applies the procedure to the known CFs mentioned 
at the beginning. In this way it is also possible to investigate how the accuracy 
depends on $K$. In this connection, we recall a theorem by 
Erde\'liy (1956) according to which, at large $q$s, the leading asymptotic term 
of  $\int_.^D (D-x)^a e^{iqx}dx$, with $a>-1$,  is $\propto e^{iqD}/q^{a+1}$ 
(confining ourselves to the only contribution related to the reported  
integration limit). It follows that the truncation order increase  makes the 
behaviour of the 3D Fourier transform of the CF approximation more accurate 
in the  region of large scattering vectors. \\  

\noindent{\bf{ 2 - Procedure for generating an approximated  CF from the CLD }}\\
\indent The mentioned mathematical structure of any polyhedron's CLD implies that 
this is analytic within each $r$-subinterval and that, within any right or left 
(small) neighbourhouds of $D_i$, its expansion reads
\begeq\label{1} 
\gamma"(r)={\gotD^{\pm}}_i(r)=\sum_{j=0}^{\infty}\Big({a^{\pm}}_{i,j}\, 
|r-D_i|^j+{b^{\pm}}_{i,j} |r-{D_i}|^{j+1/2}\Big),
\endeq
where superscript  $^+$  applies  if  $r\to {D_i}^+$ (and superscript $^-$ if  
$r\to {D_i}^-$). The above series, truncated at $j=K$,  will be denoted as 
$\gotD^{\pm}_{i,K}(r)$, \ie
\begeq\label{1a}
\gotD^{\pm}_{i,K}(r)\equiv \sum_{j=0}^{K}\Big({a^{\pm}}_{i,j}\, |r-D_i|^j+
{b^{\pm}}_{i,j} |r-{D_i}|^{j+1/2}\Big). 
\endeq
[Clearly, approximating ${\gotD^{+}}_i(r)$ by $\gotD^{+}_{i,K}(r)$ involves  
an error which is $O(|r-D_i|^{K+1})$ within a small right-neighbourhoud of $D_i$.  
The same  happens for ${\gotD^{-}}_i(r)$.]  We introduce now  the  new positive 
variables $\xi_i$ and $\eta_i$ according to the definitions 
\begeq\label{3a}
\xi_i(r)\equiv{(r-D_{i-1})^{1/2}},\quad  \eta_i(r)\equiv({D_{i}-r})^{1/2}.
\endeq 
They are not independent since they are related by 
\begeq\label{3b}
\xi_i^2+\eta_i^2 = \Delta_{i},
\endeq         
so that $\eta_i\to  \Delta_{i}^{1/2}$ if  $\xi_i\to 0$ or $r\to D_{i-1}^+$ and 
$\xi_i\to \Delta_{i}^{1/2}$ if  $\eta_i\to 0$ or $r\to {D_i}^-$.
In terms of $\xi_i$ and $\eta_i$, from \refeq{1a} follows that $\gotD^{+}_{i-1,K}(r)$ 
and $\gotD^{-}_{i,K}(r)$ respectively 
take the forms 
\begin{eqnarray}\label{3c}
\gotD^{+}_{i-1,K}(r)&=& \sum_{j=0}^{K}\Big({a^{+}}_{i-1,j}\, \xi_i^{2j}+
{b^{+}}_{i-1,j} \xi_i^{2j+1}\Big),\\
\gotD^{-}_{i,K}(r)&=& \sum_{j=0}^{K}\Big({a^{-}}_{i,j}\, \eta_i^{2j}+
{b^{-}}_{i,j} \eta_i^{2j+1}\Big)
\end{eqnarray}
so that they respectively are polynomials  of degree $(2K+1)$ of $\xi_{i}$ and $\eta_i$. 
Now, if we find a polynomial $\gotGdd_{i,K}(r)$ of $(\xi_i,\, \eta_i)$  such that 
it behaves as ${\gotD^{+}}_{i-1}(r)$  [up to the term $O(\xi_i)^{2K}$ included], as 
$r\to D_{i-1}^+$, and as ${\gotD^{-}}_{i}(r)$  [up to the term $O(\eta_i)^{2K}$ included]  as 
$r\to D_{i}^-$, according to what stated in section 1, 
$\gotGdd_{i,K}(r)$ yields an algebraic 
approximation of $\gamma"(r)$ throughout the $i$th subinterval with an error 
$\propto \Delta_i^{K}$. To determine $\gotGdd_{i,K}(r)$ we put 
\begin{eqnarray}\label{4}
 \gotGdd_{i,K}(r)\ &\equiv&\ \gotGdd_{L,i,K}(r)+\gotGdd_{R,i,K},
\end{eqnarray}
with 
\begin{eqnarray}\label{4x}
\gotGdd_{L,i,K}(r)&\equiv& \big[1+\xi^{2K+1}\,P_{L,i,K}(\eta)\big]\gotD^{+}_{i-1,K}(r),\\
\gotGdd_{R,i,K}(r)&\equiv& 
 \big[1+ \eta^{2K+1}\,P_{R,i,K}(\xi)\big]\,\gotD^{-}_{i,K}(r)\label{4y}
\end{eqnarray}
where, for notational simplicity,  we omit to append index $i$ to $\xi$ and $\eta$.  
$P_{L,i,K}(\cdot)$ and $P_{R,i,K}(\cdot)$ are two unknown polynomials to be determined. 
Contribution $\gotGdd_{L,i,K}(r)$  behaves as 
$\gotD^{+}_{i-1,K}(r)$ as $r\to {D_{i-1}}^+$ up to the term $O(\xi^{2K})$ included. 
Similarly,  as $r\to {D_i}^-$, $\gotGdd_{R,i,K}(r)$ behaves as $\gotD^{-}_{i,K}(r)$ 
up to the term $O(\eta^{2K})$. Hence, $\gotGdd_{i,K}(r)$ is the sought for 
approximation  if $\gotGdd_{L,i,K}(r)$  is $o(\eta^{2K})$ as $r\to {D_i}^-$ 
and $\gotGdd_{R,i,K}(r)$ is $o(\xi^{2K})$ as $r\to {D_{i-1}}^+$. The last 
two conditions uniquely determine the unknown polynomials. In fact, by \refeq{3b} 
and \refeq{1a},   $\gotGdd_{L,i,K}(r)$  around ${D_{i}}^-$ takes the form
\begin{eqnarray}\label{5}
&&\gotGdd_{L,i,K}(r)=\big[1+(\Delta-\eta^2)^{K}(\Delta-\eta^2)^{1/2}\,
P_{L,i,K}(\eta)\big]\times \\
&&\sum_{j=0}^{K}\Big[{a^{+}}_{i-1,2j}\, (\Delta-\eta^2)^j+
{b^{+}}_{i-1,2j+1} (\Delta-\eta^2)^j\,(\Delta-\eta^2)^{1/2}\Big]. \nonumber
\end{eqnarray}
We must require that its expansion around $\eta=0$ does not involve terms 
$O(\eta^l)$ with $l\le 2K$. This property must necessarily hold true for the expansion of 
the first factor present on the right hand side of \refeq{5}.  Since the expansions 
of $(\Delta-\eta^2)^K$ and $(\Delta-\eta^2)^{1/2}$ only involve even powers of 
$\eta$, the polynomial $P_{L,i,K}(\eta)$ must  have the form: 
$\sum_{h=0}^K p_h \eta^{2h}$. The unknown coefficients $p_0,\,p_1,\ldots,p_K$ 
are iteratively determined solving the set of equations 
\begeq\label{9}
\delta_{h,0}+\sum_{l=0}^h\frac{(-)^{h-l}p_{2l}}{\Delta^{h-l-K-1/2}}
\sum_{q=0}^{\min[K,h-l]} 
\binom{K}{q}\binom{1/2}{h-q-l}=0,\ \ h=0,1,\ldots,K,
\endeq 
resulting from the expansion of the mentioned factor.  (In the above relation, 
$\delta_{h,0}$ is the Kronecker symbol.) The remaining polynomial $P_{R,i,K}(\cdot)$ 
is determined by a similar procedure considering the expansion of $\gotGdd_{R,i,K}(r)$ 
around $\xi=0$. In this way, the algebraic function $\gotGdd_{i,K}(r)$, approximating 
the CLD within the $i$th subinterval, is fully determined.\\              
To get, within the $i$th subinterval, the corresponding algebraic approximation of the
 CF, denoted by ${\gotG}_{i,K}(r)$,  it is sufficient to integrate twice the obtained 
$\gotGdd_{i,K}(r)$,  i.e.  
\begin{eqnarray}\label{10}
{\gotG}_{i,K}(r)&\equiv&\int_{D_{i-1}}^r dx \int_{D_{i-1}}^x \gotGdd_{i,K}(y)dy + 
A_i+B_i\ r=\\
\quad &&\quad \int_r^{D_{i}} dx \int_x^{D_{i}} \gotGdd_{i,K}(y)dy + A'_i+B'_i\ r,\label{11}
\end{eqnarray}
where $ A_i,\,B_i,\,A'_i$ and $B'_i$ are arbitrary constants. We underline that the 
previous  integrals are algebraic functions because their expanded integrands only 
involve a single radical, associated either to the odd powers of $\xi$ or to the odd 
powers of $\eta$. The determination of  ${\gotG}_{i,K}(r)$ requires the 
determination of constants $A_i$ and $B_i$ or $A'_i$ and $B'_i$. This is made 
possible by the properties that the CF and its first derivative $\gamma\p(r)$ 
are continuous within the  full $r$-range $[0,\, D_M]$.  [These properties 
follow from the general integral expressions of $\gamma(r)$ and $\gamma\p(r)$, 
respectively reported by Guinier \& Fournet(1955) and  Ciccariello \etal\, (1981).]  
We ecall now a general result (Ciccariello \& Sobry, 1995) according to which the 
CLD of any polyhedron is 
a first degree $r$-polynomial in the innermost range of distances, \ie\ $[0,\,D_1]$, 
and is, therefore, fully known because the relevant constant and the slope 
respectively are  the polyhedron's angularity and  sharpness. Further, the 
angularity is related to the edges' lengths and the corresponding dihedral angles 
by equation  (4.5) of Ciccariello \etal\, (1981) , while  the sharpness is the 
sum of the contributions arising from each vertex of the considered polyhedron. 
The general expression of each of the last contributions depends on the angles 
between the edges converging into a vertex as well as on the relevant dihedral 
angles, and is given by equation (3.13) of Ciccariello and Sobry (1995).  Adding 
to these results two further properties, namely  i) $\gamma(0)=1$ and 
ii) $\gamma'(0)=-S/4V$ (related to the Porod law because $S$ and $V$  
respectively denote the surface area and the volume of  the polyhedron), we 
conclude that the CF of a polyhedron always is fully known inside $[0,\,D_1]$.   
Constants $ A_i$ and $B_i$ are determined proceeding as follows. 
From \refeq{10} and the last two mentioned properties follows that 
\begeq\label{12}
{\gotG}_{1}(r)=1-S\, r/4V+\cA  r^2/2+{\cS}r^3/3\approx \int_{0}^r dx \int_{0}^x 
\gotGdd_{1,K}(y)dy
\endeq
where $\cA$ and $\cS$ respectively denote the angularity and the sharpness. 
[We have omitted index $K$ because the first equality is exact.] 
In the same way, the continuity properties of $\gamma(r)$ and $\gamma'(r)$ at 
$r=D_M$ imply that these two functions vanish at $D_M$. Thus, from \refeq{11} 
it follows that 
\begeq\label{13}
{\gotG}_{M,K}(r)=\int_r^{D_{M}} dx \int_x^{D_{M}} \gotGdd_{M,K}(y)dy.
\endeq
Having fully determined both ${\gotG}_{1,K}(r)$ and ${\gotG}_{M,K}(r)$ we proceed 
to determining the remaining constants. Constants $ A_2$ and $B_2$, present 
in the ${\gotG}_{2,K}(r)$ definition, are uniquely determined by continuously matching 
${\gotG}_{2,K}(r)$ to ${\gotG}_{1}(r)$ at $r=D_1$, \ie\, by requiring that 
\begin{eqnarray}\label{14}
\lim_{r\to{ D_1}^+}{\gotG}_{2,K}(r)&=&\lim_{r\to {D_1}^-}{\gotG}_{1}(r) ,\\
\lim_{r\to {D_1}^+}\frac{d{\gotG}_{2,K}(r)}{dr} &=& \lim_{r\to {D_1}^-}
\frac{d{\gotG}_{1}(r)}{dr},\label{14}
\end{eqnarray}
 so that, by \refeq{10},  
\begeq\label{15}
B_2=\lim_{r\to {D_1}^-}\frac{d{\gotG}_{1}(r)}{dr} \quad {\rm and}\quad 
A_2=\lim_{r\to {D_1}^-}{\gotG}_{1}(r)-B_2\,D_1. 
\endeq 
and ${\gotG}_{2,K}(r)$ also is fully determined. Iterating the procedure, 
we successively determine ${\gotG}_{3,K}(r)$, ${\gotG}_{4,K}(r),\ldots$ and, finally, 
${\gotG}_{M,K}(r)$, making apparently useless its previous  determination reported 
in \refeq{13}. However,  each step of the recursive determination introduces an 
error and the errors sum up  as the iteration goes on. Hence, it is reasonable to 
expect that the ${\gotG}_{M,K}(r)$, obtained in the last step of the recursive chain, 
does not vanish, together with its derivative, at $r=D_M$ as it is required by 
\refeq{13}. Thus, to reduce the 
approximation errors, it is more convenient to start from ${\gotG}_{1,K}(r)$ and, 
proceeding towards the right, to successively determine  ${\gotG}_{2,K}(r),\dots,
{\gotG}_{i,K}(r)$ and then, starting from the ${\gotG}_{M,K}(r)$, given by 
\refeq{13}, and, proceeding towards the left, to successively determine  
${\gotG}_{M-1,K}(r),\ldots,{\gotG}_{i+1,K}(r)$.  For the reason already noted, it is 
extremely unlike that  ${\gotG}_{i,K}(r)$ exactly  matches ${\gotG}_{i+1,K}(r)$ at 
the point $r=D_i$. Nonetheless, their matching is still possible by suitably 
modifying the definition of one of them. To this aim, we observe that adding to 
$\gotGdd_{i+1,K}(r)$ an extra contribution of the form 
\begeq\label{16}
\gotC_{i+1,K}^{[2]}(r)\equiv (r-D_{i})^{K+1}(D_{i+1}-r)^{K+1}(\alpha+\beta\,r)=
\xi^{2(K+1)}\eta^{2(K+1)}(\alpha+\beta\,r),
\endeq
with $\alpha$ and $\beta$ arbitrary constants, the expansions of  
$[\gotGdd_{i+1,K}(r)+\gotC_{i+1,K}^{[2]}(r)]$ around $\xi=0$ or $\eta=0$ coincide 
with those of $\gotGdd_{i+1,K}(r)$ up to terms $O(\xi)^{2K+1}$ or $O(\eta)^{2K+1}$, 
respectively. Besides, 
\begeq\label{17}
\gotC_{i+1,K,L}(r)\equiv \int_r^{D_{i+1}}dx\int_y^{D_{i+1}} (y-D_{i})^{K+1}
(D_{i+1}-y)^{K+1}(\alpha+\beta\,y)
\endeq 
is such that its 2nd derivative coincides with $\gotC_{i+1,K}^{[2]}(r)$ and that it 
and its first derivative vanish at $r=D_{i+1}$.  Consequently, 
$[{\gotG}_{i+1,K}(r)+\gotC_{i+1,K,L}(r)]$ is an approximation of the CF as satisfactory 
as ${\gotG}_{i+1,K}(r)$ because it obeys all the constraints that we imposed to  
determine  ${\gotG}_{i+1,K}(r)$.  Then, to match the behaviour of ${\gotG}_{i,K}(r)$ 
at $r=D_i$, we simply substitute ${\gotG}_{i+1,K}(r)$ with 
$[{\gotG}_{i+1,K}(r)+\gotC_{i+1,K,L}(r)]$ and determine constants $\alpha$ and $\beta$, 
here present, requiring that this functions and its derivative 
respectively coincide with ${\gotG}_{i,K}(r)$ and ${\gotG}_{i,K}'(r)$ at $r=D_i$.  
Alternatively, we could modify ${\gotG}_{i,K}(r)$ instead of 
${\gotG}_{i+1,K}(r)$. To do that, we must simply add to ${\gotG}_{i,K}(r)$ the 
function 
\begeq\label{18}
{\gotC}_{i,K,R}(r)=\int_{D_{i-1}}^rdx\int_{D_{i-1}}^y (y-D_{i-1})^{K+1}
(D_{i}-y)^{K+1}(\alpha+\beta\,y)
\endeq 
and then to match $[{\gotG}_{i,K}(r)+{\gotC}_{i,K,R}(r)]$ to ${\gotG}_{i+1,K}(r)$ 
at $r=D_i$. The choice between the two possibilities depends on the values of 
${\gotG}_{i,K}(D_i^-)$ and ${\gotG}_{i+1,K}(D_i^+)$. If one of these values only is 
negative, it is the one that must be corrected because the CF cannot be negative. 
In this case then the choice is unique. In the case where both values are negative, 
the choice is dictated by the fact that the inconsistency should be as small as possible. 
Finally, in the case where both values are positive the choice presumably is that 
corresponding to a more flat behaviour of the resulting  CF approximation.  At 
this point  the explanation of a procedure able to yield an algebraic approximation 
of the CF of a polyhedron starting from the knowledge of its CLD is complete. \\
The increase of  index $K$ implies that the resulting  approximation  better  
reproduces the behaviour of the exact CLD close to each $D_i$ value  
and, simultaneously, the  3D Fourier transform of the associated CF approximation 
better reproduces, owing to the mentioned Erd\'eliy theorem, the asymptotic 
behaviour of the exact form factor  in the far asymptotic region of reciprocal 
space.  Unfortunately, as $K$ increases, the agreement improves within  intervals 
aroud the $D_i$s that generally get smaller and, in reciprocal space, the asymptotic 
behaviour sets in at larger scattering vector values. Thus, an {\em a prioir} 
estimate of the dependence of the accuracy on $K$ does not seem possible. 
 An indirect, albeit rough,  estimate  can only be obtained by analyzing the 
known CFs as reported in the following section.\\

\noindent{\bf{3 - Application to  the regular tetrahedron, octahedron and cube}}\\
\indent We have applied the described procedure to approximate the CFs of the cube, 
the octahedron and the tetrahedron stemming from their CLDs reported in the 
papers mentioned in the introduction. Figures 1 and 2 show the results obtained 
with the lowest order approximation, \ie\, $K=0$, while Fig. 3 illustrates the 
octahedron approximations for the cases $K=0,\,1$ and 2. The reader can find 
the resulting equations as well as their derivation in the deposited part. Hereafter, 
we shall confine ourselves to comment the reported figures.\\
The top panel of Fig.\,1 shows the exact and the approximated CLDs of the 
mentioned three polyhedra. In the innermost  subinterval they are exact by 
construction owing to the property derived by Ciccariello \& Sobry (1995). In the 
remaining subintervals, the approximated CLDs with $K=0$ coincide with the exact 
ones at the only end points of the subintervals. Hence, they only reproduce the 
first order discontinuities present in the exact CLD of the cube and the octahedron. 
But this property is already sufficient to reproduce the CFs with a good accuracy 
as it appears evident from the bottom panel of Fig.\,1.   
\begin{figure}[tbp]
\begin{center}
{{\includegraphics[width=8.truecm]{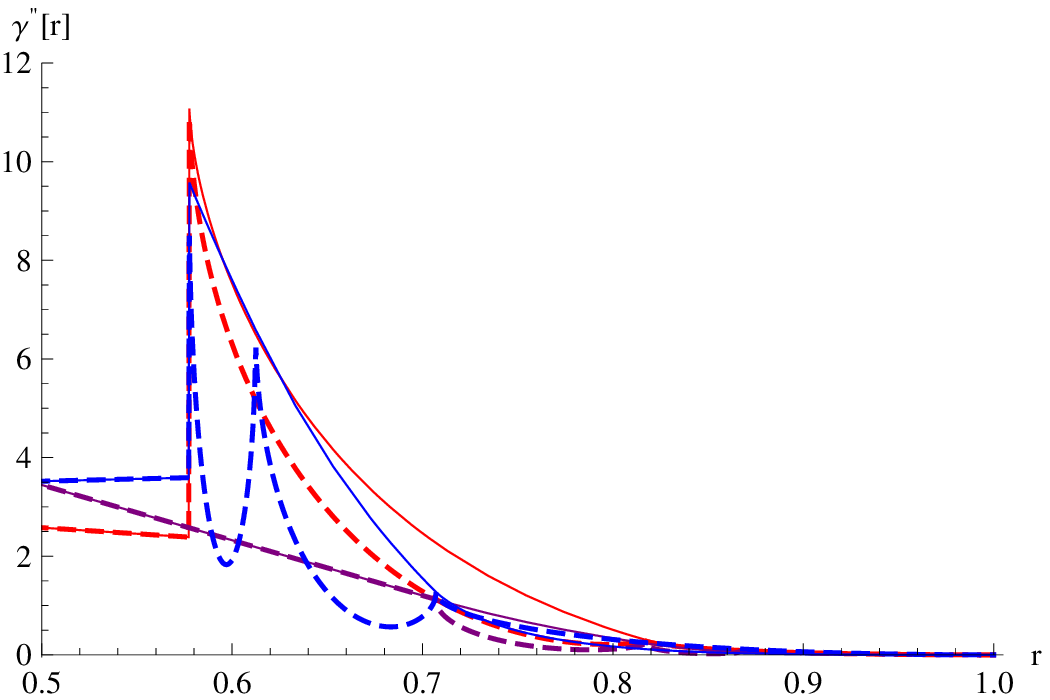}\quad\quad 
\includegraphics[width=8.truecm]{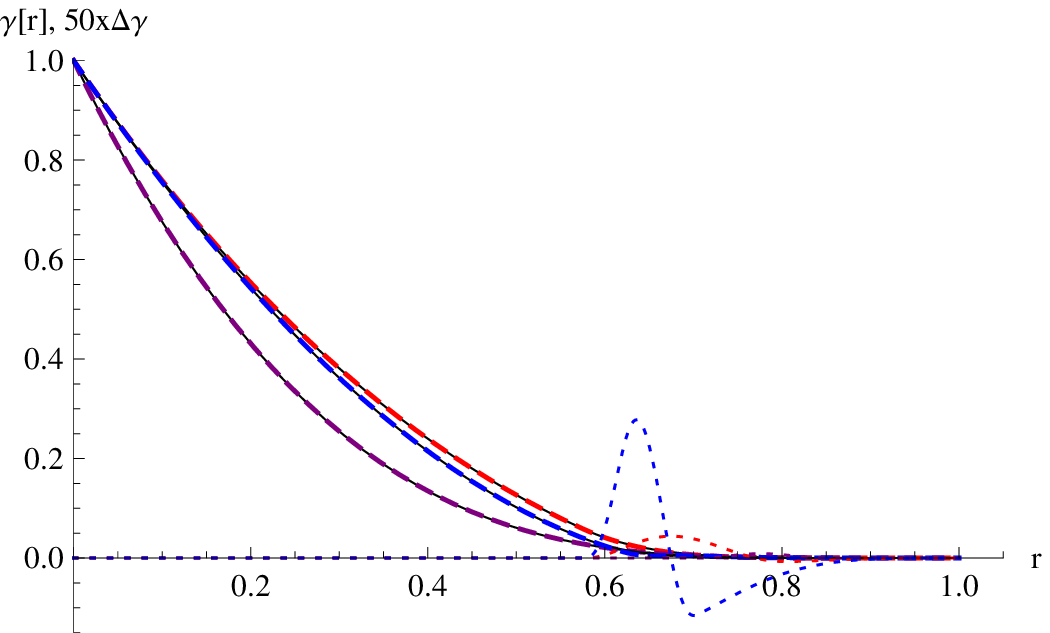}}}
\end{center}
\caption{\label{Fig1} { Top panel: the continuous curves represent the exact CLD 
of the cube (red), the octahedron (blue) and the tetrahedron (purple). The 
corresponding algebraic approximations, of order $K=0$, are represented 
by the dashed curves with the same colours.  \hfill\break\noindent 
Bottom panel: the resulting 
approximations of the associated CFs are represented with the same colours. 
The thin black lines, hardly distinguishable from the dashed ones, refer to the 
exact CFs. The dotted curves are the errors 
(multiplied by 50). }} 
\end{figure} 
This conclusion is further strenghtned by the top panel of Fig.\,2 which shows $I(q)$ 
{\em versus} $q$, \ie\, the 3D Fourier transforms (FT) of the exact and the 
approximated CFs.  
The agreement appears to be quite good throughout the reported tange of the 
scattering vector, denoted by $Q$, (instead of $q$) in the figures. We recall the 
sum-rule (Guinier\& Fournet, 1955;   Feigin \& Svergun, 1987): 
$I(0)=\int_{R^3} \gamma(r) dv= V$ where $V$ denotes the particle volume. From this 
and the fact that the approximated and the exact FTs,  for each particle shape, practically 
coincide at $q=0$, we conclude that the approximated CFs fairly  obey the  sum-rule. 
\begin{figure}
\begin{center}
{{\includegraphics[width=8.truecm]{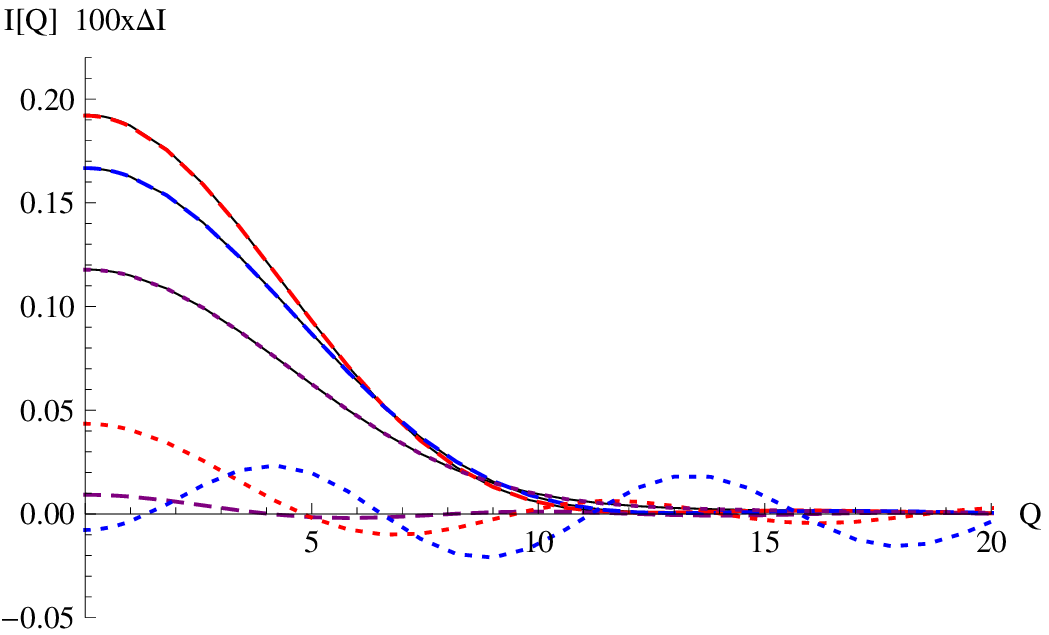}\quad\quad 
\includegraphics[width=8.truecm]{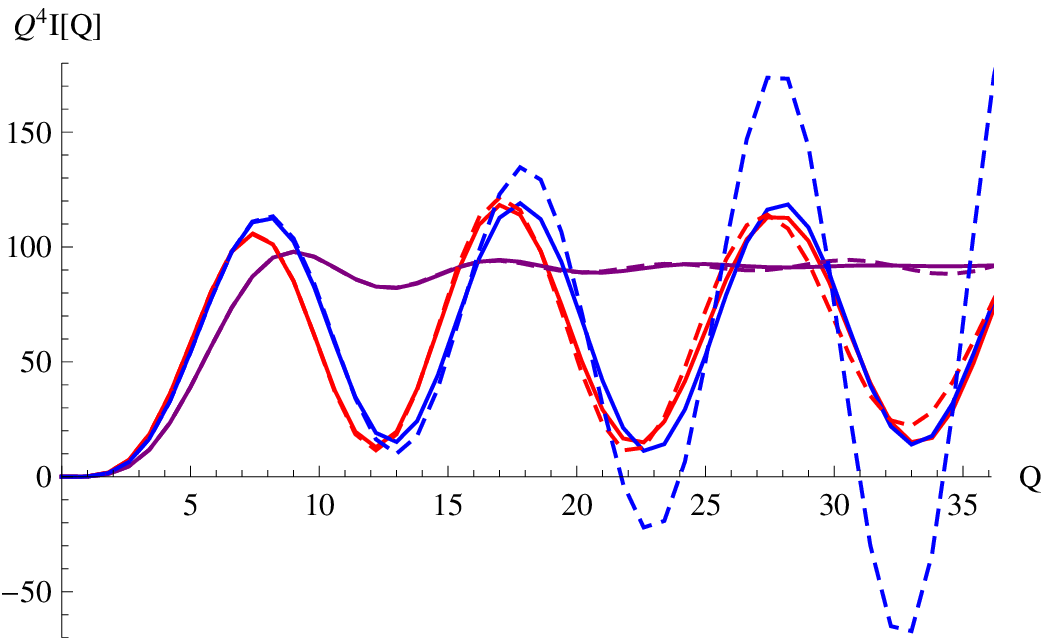}}}
\end{center}
\caption{\label{Fig2} {Top panel:  the exact and the approximated form factors of  
the three considered particle shapes, are given by the thin continuous black curves 
and the closer coloured broken ones. They practically coincide. [The colours depend on 
the particle shape as specified in the Fig.1 caption.]  
The dashed curves  are obtained by the 3D Fourier transform of  the (K=0) CF 
approximations.  The errors, multiplied by 100, are shown, with the same colours, 
by the dotted  curves oscillating around the horizonatal axis. \hfill\break\noindent
Bottom panel: Porod plot of the intensities resulting from the FTs of 
the exact (continuous) and the approximated (dashed) CFs.}}
\end{figure}    
However, the Porod plot is a tool much more accurate to evaluate the accuracy of 
an approximation. The bottom panel shows the Porod plots of the considered 
approximations. One sees that the scattering intensities relevant to the CF 
approximations of the tetrahedron and the cube are accurate throughout the 
reported $q$-range, while that of the octahedron is only accurate up to 
$q\approx 10$. In the three cases, however, one notes that the accuracy deteriorates 
as $q$ increases. This is by no way surprising because the $K=0$ approximations of 
the CLDs only reproduces the first-order discontinuities (Ciccariello, 1985) of the exact CLDs. The CLDs' 
higher order derivatives show further singularities [see Ciccariello (205b, 2014b)] that 
are responsible for further damped oscillatory contributions in the Porod plots. 
Besides, the CLD approximations, shown in Fig.\,1, show artificial wells that, by 
Fourier transforming, yield peaks around $q=2\pi/\delta$ with $\delta$ equal to 
the positions of the minima of the wells. Only beyond the largest of these $q$ values, 
the approximated and the exact Porod plots are expected to coincide and 
the noted discrepancies to disappear. \\  
\begin{figure}
\begin{center}
{{\includegraphics[width=8.truecm]{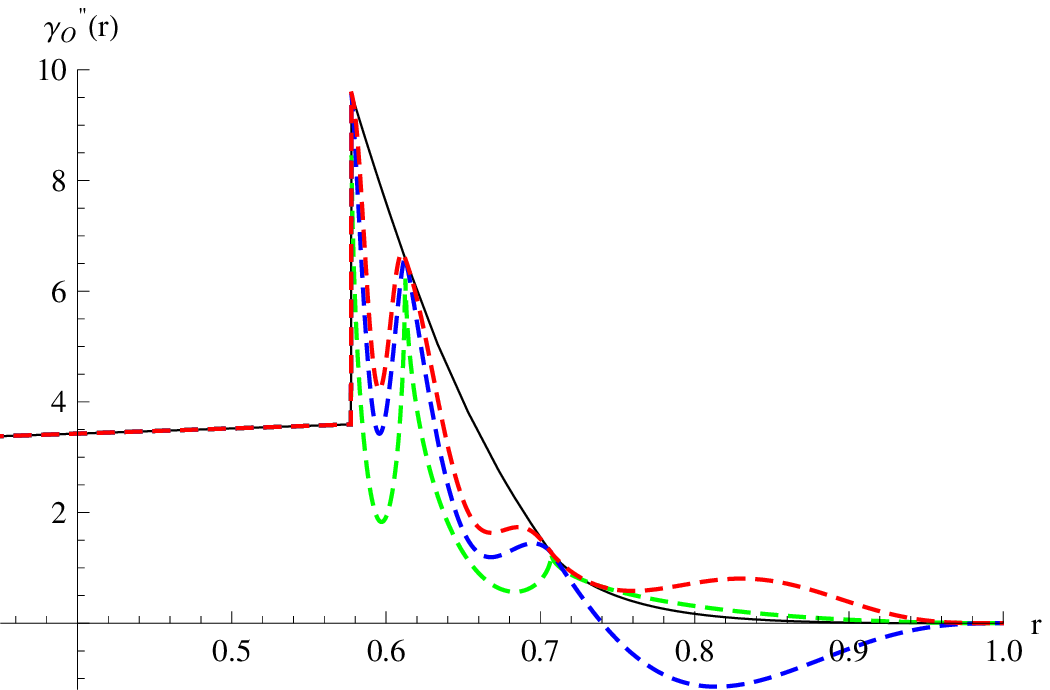} \quad\quad 
\includegraphics[width=8.truecm]{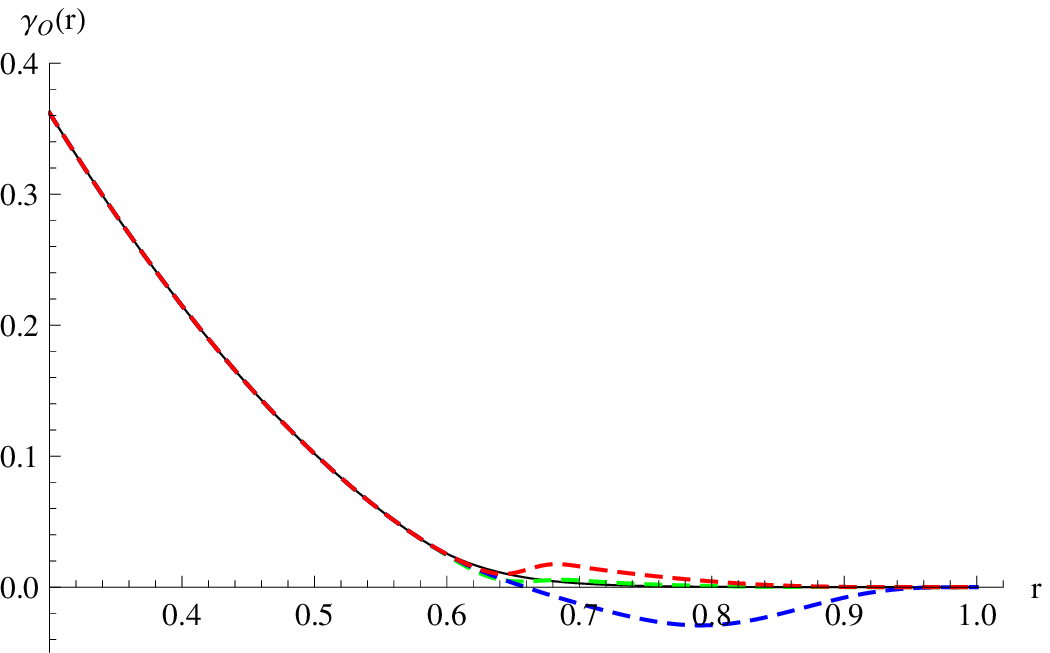}\quad\quad
\includegraphics[width=8.truecm]{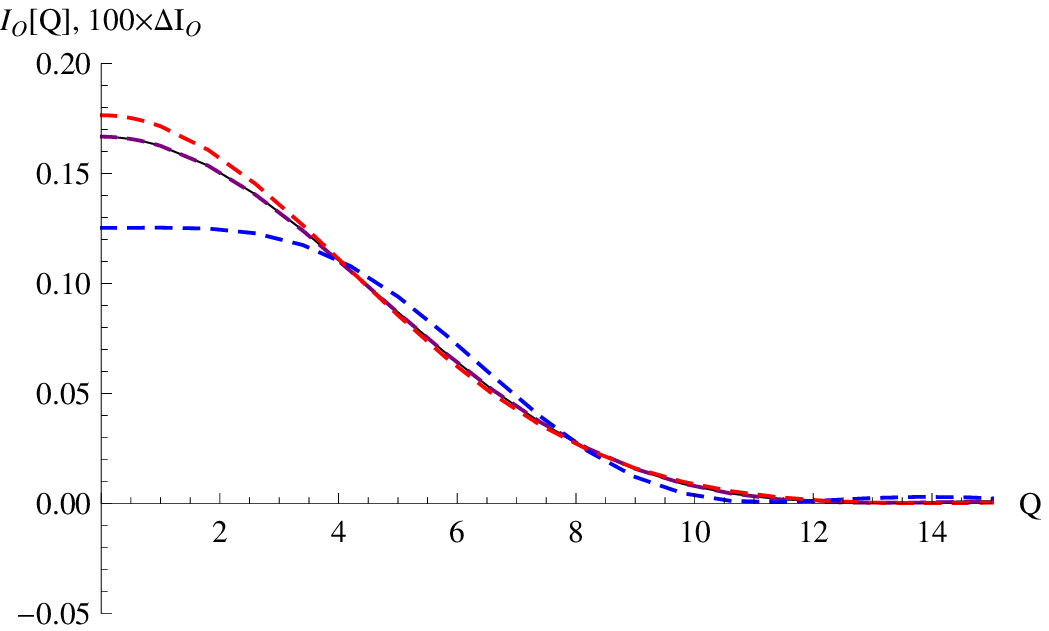}}}
\end{center}
\caption{\label{Fig3} {Top panel: the thin black curve is the exact CLD and the broken  
green, blue and red curves respectively represent the CLD approximations relevant 
to $K=0$, $K=1$ and $K=2$. \hfill\break\noindent  
 Middle panel: with the same colours we plot the CF approximations resulting from 
the above CLD approximations.  \hfill\break\noindent
Bottom panel: 
plots of the scattering intensities relevant to the three approximations in the inner 
$q$-range. }}
\end{figure} 
Fig.\,3 allows us to appreciate how things change as we increase the approximation 
order. It refers to the only octahedron which has a CLD more structured than the 
tetrahedron's and the cube's. The top panel shows that, as $K$ increases from 0 to 2, 
the approximated CLD becomes, so to speak, more adherent to the exact one 
around the end points of the distance subintervals. In the two internal subintervals, 
the approximated CLD gets nearer to the exact one throughout the full subintervals, 
while in outer one it gets farther as we pass from $K=0$ to $K=1$ and then closer 
for $K=2$, remaining however farther than the $K=0$ approximation. This last 
discrepancy propagates towards the inner two subintervals owing to the matching 
procedure so that the final accuracy of the total CF approximation worsens as we 
pass from $K=0$ to $K=2$ and to $K=1$, as it appears in the middle panel. This 
is confirmed by the bottom panel that shows the corresponding scattering intensity 
in the innermost $q$-range. The improvement of the approximations as $K$ increases 
can only be appreciated by the corresponding Porod plots that are reported in  the 
ms'  part deposited with IUCr. There it appears that  the $K=1$ and $K=2$ intensities almost 
coincide with the exact one in the region $q>500$, while  we must go beyond 
$q=3000$ for this to happen for the $K=0$ approximation.\\ 

\noindent{\bf{4 - Conclusions}}\\
\indent From the above results it appears reasonable to conclude that the simplest 
approximation, relevant to the choice $K=0$, yields an algebraic approximation 
of the CF accurate enough to meet the standards of crystallographers and 
small-angle scattering  people. We conclude with two remarks. First, the reported 
procedure still works if the  value of $K$  is  differently chosen in the different 
subintervals. For instance, the behaviour of the CLD approximations, shown 
in Fig.\,2, suggests that the choice $K=2$ in the second and third subinterval 
and $K=0$ in the fourth ought to be more accurate because the resulting 
approximation is closer 
to the exact CLD. Second, in constructing the CF approximation, the crucial 
point is that the approximation must continuously interpolates, together with 
its derivatives,  the truncated expansions 
of the CLD at the end-points of the the considered subinterval. We have 
illustrated a procedure that achieves the aim, but other procedures are possible. 
For instance, taking advantage of the fact that the approximation error reduces 
with the subinterval lenghts, one could divide each subinterval into three parts, 
approximate the CLD with its  left and right truncated expansions in the first and 
the third of these intervals, then continuously interpolate the truncated expansions, 
evaluated at the dividing points, by the above explained procedure and, finally, 
get the full CF algebraic approximations by the matching procedure. The paid 
cost is the greater complexity of the approximation.  \\ 

\noindent{\bf{Acknowledgments}}\\
\indent I gratefully thank the anonymous referee  for having raised
the question of the possibility of approximating the CF starting from the CLD
knowledge. 
\vfill\eject
\subsection*{References}
\begin{description}
\item[\refup{}]  Ciccariello, S. (1985).  {\em  Acta Cryst.}  A{\bf 41}, 560-568.
\item[\refup{}]  Ciccariello, S. (2005a).  {\em J. Appl. Cryst.} {\bf 38}, 97-106.
\item[\refup{}] Ciccariello, S.  (2005b). {\em Fibres \& Text. East Eur.} {\bf 13}, 41-46.
\item[\refup{}] Ciccariello, S. (2014a). {\em J. Appl. Cryst.} {\bf 47}, 1216-1227.
\item[\refup{}] Ciccariello, S. (2014b). {\em J. Appl. Cryst.} {\bf 47}, 1445-1448.
\item[\refup{}] Ciccariello, S. (2020a). {\em arXiv}:1911.02532v2 [math-ph].
\item[\refup{}] Ciccariello, S. (2020b). {\em Acta Cryst.} A{\bf 76}, 
https://doi.org/10.1107/S2053273320004519.
\item[\refup{}] Ciccariello, S.,   Cocco, G.,  Benedetti, A. \&   Enzo, S. (1981).
{\em  Phys. Rev.} {\bf B23}, 6474-6485.  
\item[\refup{}] Ciccariello, S. \& Sobry, R.  (1995). {\em Acta. Cryst.} {\bf A51}, 60-69.
\item[\refup{}] Erd\'eliy, A. (1956). {\em Asymptotic Expansions}, ch. II. New York:Dover.
\item[\refup{}] Feigin, L.A.  \&   Svergun, D.I.  (1987). {\em Structure Analysis
by Small-Angle X-Ray and  Neutron Scattering}, New York: Plenum Press.  
\item[\refup{}] Goodisman,  J. (1980). {\em  J. Appl. Cryst.} {\bf 13}, 132-134.
\item[\refup{}] Guinier, A.  \&  Fournet, G.   (1955). {\em Sall-Angle Scattering of 
X-rays}, New York: John Wiley.
\end{description}
\end{document}